\documentclass[a4paper,12pt,reqno]{amsart}
\usepackage{amsmath}%
\usepackage{amsfonts}%
\usepackage{amssymb}%
\usepackage{graphicx}

\theoremstyle{plain}

\numberwithin{equation}{section}

\begin{document}
\title[Veracidad de CH. Versi\'on en espa$\tilde{N}$ol]{Algunas ideas acerca de la numerabilidad de los conjuntos infinitos}
\author{Denis Mart\'inez T\'apanes}
\email[A. One]{denismt@ucm.vcl.sld.cu}%
\thanks{This paper is in final form and no version of it will be submitted forpublication elsewhere.}
\date{Diciembre 31,2015}
\dedicatory{Dedicado a mi profesor Lorgio Batard.}

\maketitle
\begin{quotation}
\small R\tiny{ESUMEN}
\footnotesize{Donde se hace una revisi\'on del concepto de numerabilidad en las matem\'aticas, sometiendo a cr\'itica algunos de los teoremas hasta ahora aceptados, demostrando su inconsistencia y dando, adem\'as, elementos concretos sobre la numerabilidad de todas las potencias del conjunto de los n\'umeros naturales.}\\

\small A\tiny{BSTRACT}
\footnotesize{In which a review of the concept of countability is done in mathematics, subjecting review some of the theorems so far accepted, showing their inconsistency and also taking concrete elements on the countability of all the powers of the set of natural numbers.}
\end{quotation}
\section{Introduction}
Como bien es conocido en el transcurso de todo el siglo XX se produjeron las tentativas de demostrar la hip\'otesis del cont\'inuo. Prime-\;\;\;\;\;\;ramente se supuso que cantor estaba  equivocado, luego se vi\'o que la hip\'otesis del continuo era falsa, pero m\'as tarde se descubri\'o  que la demostraci\'on de Julius K\"onig  conten\'iía un error (lo descubri\'o Zermelo).\\En $1940$ Kurt G\"odel demostr\'o que la hip\'otesis del continuo no puede ser refutada en el marco del sistema axiom\'atico $ZFC$. Para ello, G\"odel a$\tilde{n}$adi\'o la propia hip\'otesis del continuo como axioma a los de $ZFC$ y demostr\'o que se obten\'ia un sistema consistente. Por otra parte, Paul Cohen demostr\'o en $1963$ que la hip\'otesis del continuo no puede ser demostrada en $ZFC$ a$\tilde{n}$adiendo el contrario de la hip\'otesis del continuo a $ZFC$ y demostrando, como G\"odel, que el sistema de axiomas que se obten\'ia era de nuevo consistente. Es decir, la hip\'otesis del continuo es independiente de $ZFC$, lo que significa que se puede construir una teor\'ia de conjuntos consistente donde la hip\'otesis del continuo sea cierta y tambi\'en puede construirse una teor\'ia de conjuntos consistente donde dicho resultado sea falso.\\Este trabajo se  conduce por otra v\'ia muy diferente a la que llev\'o a los matem\'aticos del pasado siglo al callej\'on sin salida de la indecidibilidad de la famosa hip\'otesis de Cantor.\\ 
En fin de cuentas $ZFC$ mostr\'o ser sumamente inseguro para llegar a la verdad de la teor\'ia, puesto que pod\'ian ser ``encajadas'' en su seno afirmaciones que son negaci\'on una de la otra manteni\'endose ellos, en un caso y en otro, perfectamente funcionales para la teor\'ia de conjuntos. Tal opini\'on era compartida por G\"odel y Cohen aunque este \'ultimo era m\'as formalista. A todas luces lo \'unico a que se ha llegado es a que la hip\'otesis del cont\'inuo es, por igual, verdadera y falsa, lo cual no se puede considerar satisfactorio como soluci\'on cient\'ifica. Incluso existe la opini\'on de que dif\'icilmente, Hilbert hubiera aceptado tales conclusiones como respuesta a su primer problema.\\ 
Por ello, en vez de plantear un sistema de axiomas y partir de \'el, demuestro directamente que todo conjunto tiene la misma potencia que el conjunto de todos sus subconjuntos. Es decir, en vez de analizar la hip\'otesis del continuo suponiendo que ella se desprende l\'ogicamente de una teor\'ia en la que no hay errores, demuestro que precisamente los errores se produjeron antes de arribar a la famosa hip\'otesis.\\
La ``verdad'' acerca de la diferencia de cardinalidad de un conjunto y su conjunto potencia no era una tal verdad y de all\'i por supuesto toda la confusi\'on posterior. Desde mi punto de vista es natural que prime-\\ramente se concluyera que no exist\'ian potencias intermedias entre $\mathbb{N}$ y $\mathbb{R}$ puesto que el teorema 2, aqu\'i demostrado, muestra la numerabilidad de todas las potencias a partir de $\mathbb{N}$. \\
En la segunda secci\'on de este trabajo se muestra con claridad varios de los errores que conducen a la creencia de que el conjunto de los reales tiene una potencia mayor que el de los naturales, incluyendo la inconsistencia del m\'etodo diagonal de Cantor aplicado con este fin. Los errores l\'ogicos que aqu\'i se se$\tilde{n}$alan forman parte de la exposici\'on com\'un de la teor\'ia de conjuntos y por tanto se han propagado a la l\'ogica aunque es evidente que su impacto no ha llegado a las matem\'aticas verdaderamente importantes por su aplicabilidad que han demostrado, en la pr\'actica, su val\'ia como herramienta de las ciencias, sino que se ha quedado en el campo de la especulaci\'on filos\'ofica de los neo-positivistas.\\
Curiosamente las matem\'aticas nunca estuvieron tan ``inseguras'' hasta el surgimiento de la rama que deb\'ia darles la mayor seguridad, es decir la Teor\'ia de Conjuntos, que constituye supuestamente su base te\'orica general, su fundamento. Desde mi punto de vista, si en verdad la teor\'ia creada por Cantor fuera el fundamento de todas las matem\'aticas, y no solamente la teor\'ia que los matem\'aticos DESEAN que sea dicha base, toda la matem\'atica no ser\'ia m\'as que un c\'umulo inmenso de afirmaciones inservibles; para darse cuenta de ello solamente hace falta mirar el estado de cosas presente en ese campo lleno de desacuerdos, preguntas sin respuesta, y respuestas que no satisfacen verdaderamente a nadie como es el caso del tema relacionado con la hip\'otesis del continuo.\\   
En la secci\'on $3$ proponemos los teoremas que resuelven el estado contradictorio de cosas que revela la secci\'on $2$; y, por \'ultimo se dan las conclusiones iniciales al menos sobre los resultados obtenidos.
\section{An\'alisis cr\'itico de algunos teoremas cl\'asicos}
\subsection{Primer caso}
La primera demostraci\'on a que queremos hacer referencia es la dada en el teorema $1$ del ep\'igrafe $2$, de [1] sobre la innumerabilidad de los reales.\\
En dicha demostraci\'on se realiza expl\'icitamente la suposici\'on de la existencia de una lista que contiene todos, o una parte, de los n\'umeros reales del segmento $[0,1]$ y luego se supone un n\'umero construido de una manera determinada\footnote{M\'etodo diagonal de Cantor.} que se diferencia de todos los elementos de la lista dada.\\
La construcci\'on de tal n\'umero se realiza diferenci\'andolo sucesivamente de los primeros elementos de la lista pero no tiene en cuenta que la lista, al contener, por hip\'otesis, todos\footnote{Por otra parte, si en la lista no est\'an todos los del segmento $[0,1]$ es evidente que se puede construir uno de este segmento que no est\'a en la lista.} los elementos posibles, siempre tendr\'a un n\'umero que tiene sus primeras cifras decimales iguales al n\'umero que se construye y que de la misma manera que el n\'umero en construcci\'on se diferencia de todos los anteriores, ser\'a igual, hasta donde ha sido construido, a una cantidad infinita de n\'umeros posteriores; de manera que el proceso constructivo nunca puede llegar a obtener un n\'umero que ya no est\'e "m\'as abajo".\\
El principio de inducci\'on transfinita, incluso, puede ayudar a demostrar la imposibilidad de construir un n\'umero que no est\'a en la lista dada: El proceso de construcci\'on genera una serie infinita de n\'umeros racionales  $S={(0,b_1 ),(0,b_1 b_2 ),\ldots,(0,b_1 b_2…b_n )}$, pero si el n\'umero $0,b_1 b_2…b_n$ est\'a en la lista (lo cual es indudable) estar\'a evidentemente $0,b_1 b_2…b_n b_{(n+1)}$, por lo que puede asegurarse que todos los  elementos de la sucesi\'on $S$ est\'an en la lista incluyendo su punto l\'imite.\\
Otro aspecto interesante es que la demostraci\'on de este teorema en realidad usa dos condiciones independientes, por lo que tiene la forma $A\wedge B\Rightarrow C$, donde $A=$(Existencia de la lista), $B=$(Existencia de un n\'umero construido de una forma dada), y $C=$(Contradicci\'on). Por lo que lo correcto es lo contrario a la conjunci\'on de suposiciones que se ha hecho, es decir que lo cierto es
\begin{center}
$\neg A\vee \neg B=$ (No existe la lista, o de lo contrario\\ no puede construirse un n\'umero de la forma dada)
\end{center}
, pero lo anterior es cierto, como es sabido, si se cumple cualquiera de las alternativas, por lo que no demuestra realmente $\neg A$.\\ 
De esto se desprende que en la demostraci\'on de cierto $A$ por reducci\'on al absurdo, se supone $\neg A$ y debe realizarse un razonamiento de la forma,
\[(\neg A\Rightarrow B_1 )\wedge(B_1\Rightarrow B_2 )\wedge(B_2\Rightarrow B_3 )\wedge \ldots\wedge (B_n\Rightarrow Absurdo)\]
, lo cual es equivalente a
\[(A\Leftarrow\neg B_1 )\wedge(\neg B_1\Leftarrow\neg B_2 )\wedge\ldots\wedge(\neg B_{n-2}\Leftarrow\neg B_{n-1})\wedge(\neg B_n\Leftarrow No\, Absurdo)\]
\subsection{Segundo caso}
La segunda demostraci\'on que queremos analizar es la dada en el teorema $1$ del ep\'igrafe $3$ cap\'itulo $V$ de $[2]$ sobre la innumerabilidad de los reales.\\ 
Aqu\'i es necesario aclarar que la aplicaci\'on de los l\'imites no es v\'alida en este caso, debido a que, una desigualdad como la dada all\'i $c\neq a_n$ puede cumplirse para todo valor finito de $n$ pero eso no implica que $c\neq{lim}_{n\rightarrow\infty}(a_n) $ (N\'otese que $1\neq1+1/n,n<\infty$, sin embargo $1={lim}_{n\rightarrow\infty}(1+1/n))$.\\
Pero no s\'olo eso, sino que el m\'etodo de construcci\'on del n\'umero $c$ es tal que los intervalos elegidos sucesivamente, y cuya intersecci\'on es precisamente $c$, son intervalos que necesariamente pertenecen a $[0,1]$; por tanto si bien es cierto que $a_{n}\notin p_{n}q_{n}$ tambi\'en es cierto que todos los elementos de $p_{n}q_{n}$ pertenecen a $[0,1]$, ya que por definici\'on todos los elementos de $[0,1]$ est\'an en la lista (esto, por supuesto para todo valor de $n$), y por consiguiente el punto intersecci\'on de todos estos sub-intervalos pertenece, tambi\'en, a la lista. De manera que este procedimiento conduce a un punto que en realidad est\'a en la lista y no fuera de ella como se quiere mostrar.
\subsection{Tercer caso}
Con respecto a la afirmaci\'on ampliamente aceptada de que $c=2^{\aleph_0}$, podemos hacer referencia a la demostraci\'on presentada en el ep\'igrafe $4$ del capítulo $VI$, de $[2]$. All\'i se supone que se ha establecido una correspondencia $uno-a-uno$ entre los conjuntos de puntos reales del intervalo $(0,1)$ y el conjunto de todas las secuencias infinitas $t=(t_1,t_2,…)$, conformadas por ceros y unos, mediante la funci\'on
\[f(t=(t_1,t_2,…))=\left\{
\begin{matrix}
\sum_{i=1}^\infty (t_i/2^i)           ,t \,presenta\, una\, cantidad \,infinita \,de \,ceros \\
1+\sum_{i=1}^\infty (t_i/2^i)         ,t \,presenta \,una \,cantidad \,finita \,de \,ceros
\end{matrix}
\right\}\]
N\'otese, sin embargo, que el n\'umero $0,3$ es representante del conjunto infinito de los elementos de $(0,1)$ que no puede estar en la imagen de esta funci\'on, en virtud de que $0,3=3/(5*2)$ mientras que, ninguna de estas sumas puede tener, por definici\'on, un denominador que no sea, exclusivamente, potencia de dos (Por ejemplo $0,125$ es admisible). As\'i queda claro que, al menos, esta demostraci\'on no es v\'alida, tal y como las que hemos analizado hasta el momento.
\subsection{Cuarto caso}
Por \'ultimo, mostremos la inconsistencia de la demostraci\'on sobre la imposibilidad de establecer una correspondencia biun\'ivoca entre los elementos de cierto conjunto y los de su conjunto potencia correspondiente\footnote{Constitye una versi\'on del teorema de Cantor.}, que aparece en $[1]$. La misma se basa en la existencia de la siguiente situaci\'on formal dado un conjunto $B$ arbitrario:
\[\exists A\subset B,\forall D\subset B,\forall x\in B,((x\leftrightarrow D)\wedge(x\notin D)\Rightarrow  \]
\begin{equation}
\Rightarrow(x\in A))\wedge((x\leftrightarrow D)\wedge(x\in D)\Rightarrow(x\notin A))   
\end{equation} 
Donde $x\leftrightarrow D$ es cierta correspondencia que se ha establecido entre los elementos de $B$ y sus subconjuntos. El hecho es que cuando deseamos conocer el estatus del propio conjunto $A$ debemos sustituir $D$ por $A$ y se obtiene a partir de $(2.1)$ lo siguiente:
\[\exists A\subset B,\forall A\subset B,\forall x \in B,((x \leftrightarrow A)\wedge(x \notin A)\Rightarrow\]
\[\Rightarrow(x \in A))\wedge((x \leftrightarrow A)\wedge(x \in A)\Rightarrow(x \notin A))\] 
O sea, se ha obtenido una situaci\'on parad\'ogica que nos lleva a concluir que no tiene lugar que $\exists x\in B,(x\leftrightarrow A)$ y por tanto no se puede establecer una correspondencia entre los elementos de $B$ y todos sus subconjuntos.\\ 
Veamos el asunto, ahora, un poco m\'as profundamente: Si damos por sentada, inicialmente, la existencia de la relaci\'on $uno-a-uno$ entre los elementos de $B$ y sus subconjuntos, ser\'a que $(x\leftrightarrow D)\wedge(x\notin D)\equiv(x\notin D)$ y $(x\leftrightarrow D)\wedge(x\in D)\equiv(x\in D)$, por lo que $(2.1)$ quedar\'a de la siguiente forma
\[\exists A \subset B,\forall D \subset B,\forall x \in B,((x \notin D)\Rightarrow(x \in A))\wedge((x \in D)\Rightarrow(x \notin A))\]
O, lo que es lo mismo
\[\exists A\subset B,\forall D\subset B,\forall x\in B,((x\notin D)\Leftrightarrow(x\in A))\]       
Pero ahora se revela algo fundamental que tiene lugar, y es que, de lo anterior, se deduce, necesariamente, que $A\cap D=\emptyset$. Esto es consecuencia de nuestras suposiciones iniciales y no puede ser violado, por lo que es totalmente ilegal, desde el punto de vista l\'ogico, el cambio de $D$ por $A$ que produce la contradicci\'on en virtud de la cual se demuestra el teorema. De hecho lo que en verdad se tiene es que
\[((x\notin D)\Leftrightarrow(x\in A))\Rightarrow (A\cap D=\emptyset)\Rightarrow (A\neq D)\]    
Por lo que lo que la igualdad de $A$ y $B$ conduce necesariamente a lo que sigue
\[((x\in D)\Leftrightarrow(x\in A))\Leftarrow (A\cap D\neq\emptyset)\Leftarrow (A= D)\]
Es decir, no se produce la paradoja, y para ello solamente es necesario respetar la propia l\'ogica formal\footnote{En el teorema $2$ de la siguiente secci\'on se demuestra de manera contundente que el teorema de Cantor (Cardinal de x$\neq$ Cardinal de P(x) ) es completamente falso. }.
\section{Propuestas te\'oricas concretas ante la situaci\'on planteada en la secci\'on $2$}
\textbf{\emph{Definici\'on.}}\textit{Un $P-conjunto$ $(Pc)$ es un conjunto que posee una propiedad $P$ dada que garantiza que:1-Un conjunto $Pc$ puede representarse como la uni\'on de cualquier cantidad, finita o numerable, de otros conjuntos que tambi\'en son $Pc$ y toda uni\'on de cualquier cantidad, finita o numerable, de conjuntos $Pc$ es también $Pc$.
2-Si un conjunto $A$ es no $Pc$ entonces el conjunto $A-B$, donde $B$ es finito y $B\subset A$, es no $Pc$.}\\
\textbf{Teorema 1.}\textit{Existe al menos un $Pc$ infinito que no tiene ning\'un subconjunto $Pc$ diferente de \'el mismo.}\\
\textbf{Demostraci\'on.}Supongamos lo contrario. Entonces un conjunto infinito $A$, arbitrario pero que no es $Pc$, tendr\'a un subconjunto $Pc$ y podr\'a por tanto representarse como
\[A=Pc_1\cup A_1\]
El conjunto $A_1$ es infinito, ya que de lo contrario se pudiera escribir $A-A_1=Pc_1$, donde el miembro derecho es $Pc$ mientras que el izquierdo, en virtud de lo definido en el punto $3$, no es $Pc$ y se generar\'ia una contradicci\'on.\\
Pero  $A_1$ debe tener tambi\'en un subconjunto $Pc$ diferente de \'el mismo como consecuencia de nuestro supuesto, por lo que podemos escribir
\[A=Pc_1\cup Pc_2 \cup A_2\]
Donde $A_2$ es infinito, ya que de lo contrario se pudiera escribir $A-A_2=Pc_1\cup Pc_2$, donde el miembro derecho es $Pc$, en virtud de lo planteado en el punto $2$, mientras que el izquierdo, gracias a lo definido en el punto $3$, no es $Pc$ y se generar\'ia una contradicci\'on.\\
El procedimiento puede repetirse $n$ veces manteniendo las mismas suposiciones y obtenerse
\[A=(\bigcup_{i=1}^n Pc_i )\cup A_n\]
En el l\'imite, cuando $n\rightarrow\infty$, se suponen ``extra\'idos'' todos los $Pc$ pero la situaci\'on debe mantenerse quedando una representaci\'on an\'aloga a la anterior, es decir
\[A=(\bigcup_{i=1}^\infty Pc_i )\cup(lim_{n\rightarrow\infty}A_n)\]
Donde ($lim_{n\rightarrow\infty}A_n$ ) no puede ser un conjunto finito, ya que de lo contrario se pudiera escribir $A-(lim_{n\rightarrow\infty}A_n )=(\bigcup_{i=1}^\infty Pc_i )$, donde el miembro derecho es $Pc$, en virtud de lo planteado en el punto $2$, mientras que el izquierdo, gracias a lo definido en el punto $3$,  es no $Pc$ y se generar\'ia como en el caso anterior, una contradicci\'on.\\
De aqu\'i se deduce que ($lim_{n\rightarrow\infty}A_n$) es un conjunto infinito y no posee un subconjunto $Pc$ ya que el procedimiento de ``extracci\'on'' de subconjuntos $Pc$ se ha llevado al l\'imite en que todos pasaron a la uni\'on ($\bigcup_{i=1}^\infty Pc_i$ ). Sin embargo esto contradice nuestra suposici\'on de que todo conjunto infinito tiene un subconjunto $Pc$. El teorema est\'a demostrado.\\
En la parte correspondiente a las conclusiones se expondr\'an las consecuencias profundas de este teorema.\\
En relaci\'on con la numerabilidad de los reales proponemos, ahora, lo siguiente:\\
\textbf{ Teorema 2.} \textit{Si cierto conjunto $A$ es numerable, entonces su conjunto potencia $P(A)$ es igualmente numerable.}\footnote{Esto constituye un contraejemplo al teorema de Cantor que demuestra su falsedad.}\\
\textbf{Demostraci\'on.} Se tiene que, por ser $A$ numerable, se puede poner en forma de lista
\[A=\left\{\alpha_1, \alpha_2, \alpha_3,\ldots, \alpha_n,…\right\}\]
Entonces el conjunto $P(A)$ se puede ordenar de la siguiente forma: Definamos la altura $h$ del subconjunto dado como la suma de los sub\'indices de sus elementos; es decir, la altura de $\left\{\alpha_1\right\}$ es igual a la unidad ``$1$'', mientras que la altura de $\left\{\alpha_1,\alpha_2\right\}$ es igual a ``$3$''. Para cada altura en particular se tendr\'a una secuencia evidentemente finita $S_h$ de subconjuntos de la misma altura; por ejemplo ser\'a
\[S_1=\left\{\left\{\alpha_1\right\}\right\}\]
\[S_2=\left\{\left\{\alpha_2\right\}\right\}\]
\[S_3=\left\{\left\{\alpha_3\right\},\left\{\alpha_1,\alpha_2\right\}\right\}\]
\[\vdots\]
\[S_{20}=\left\{\left\{\alpha_{20}\right\} ,\left\{\alpha_1,\alpha_{19} \right\},\ldots,\left\{\alpha_1,\alpha_2,\alpha_{17} \right\},\ldots,\left\{\alpha_1,\alpha_2,\alpha_3,\alpha_{14} \right\},\ldots,
\left\{\alpha_2,\alpha_3,\alpha_4,\alpha_5,\alpha_6\right\}\right\} \]
\[\vdots\]
Ahora solo falta poner $A$ como la lista de sucesiones seg\'un la altura
\[P(A)=\left\{\left\{\emptyset\right\},S_1,S_2,\ldots,S_n,\ldots\right\}\] 
El teorema est\'a demostrado.\\ 
De m\'as est\'a decir que, si este teorema es incorrecto (y es bastante evidente que lo que ocurre es lo contrario), la numerabilidad de los racionales queda anulada debido a que el esquema l\'ogico de ambas demostraciones es sumamente an\'alogo. 
\subsection{M\'etodo de conteo para los reales}
Por \'ultimo, vamos a proponer una manera de ``contar'' los reales del intervalo $\left[0,1\right]$, en lo que hemos denominado ``Procedimiento de Conteo por Capas'' (Layered count procedures): La primera capa est\'a conformada por $10$ elementos 
\[0,0\ldots 0,1\ldots 0,2\ldots\ldots\ldots 0,9\]
La segunda capa est\'a constituida por $100$ elementos generados introduciendo en cada uno de los anteriores los d\'igitos desde el cero hasta el  nueve en la segunda cifra despu\'es de la coma. Por ejemplo ser\'a
\[\begin{matrix}
0,0&0,1\ldots&\ldots0,9 \\
(0,00\;\;\;\;0,01\;\ldots\;0,09)&(0,10\;\;\;\;0,11\;\ldots\;0,19)\ldots&(0,90\;\;\;\;0,91\;\ldots\;0,99)
\end{matrix}\]
As\'i se puede continuar asignando un n\'umero natural a cada elemento de cada capa de izquierda a derecha y comenzando por la primera capa. Est\'a claro que cada  elemento  se autogenera en la capa posterior por lo que en la capa de altura $h$ solo ser\'a necesario numerar $10^h-10^{h-1}=9*10^{h-1}$ elementos a exepci\'on de la primera. Es indudable que en esta lista, por capas, estar\'an todos reales de $\left[0,1\right]$ sin que falte ninguno de ellos a pesar de que los irracionales estar\'an confinados en la capa correspondiente al paso al l\'imite $ n\longrightarrow\infty$; pero los naturales siempre ser\'an suficientes no importa la altura de la capa.
\section{Conclusiones}
El problema de los $P- conjuntos$ tiene una importante consecuencia para los fundamentos de las matem\'aticas debido a que, si se cambia la denominaci\'on de $P- conjuntos$  por ``\emph{conjuntos numerables}'' y la propiedad $P$ se identifica con la igualdad con el conjunto de los naturales en cuanto a potencia se refiere, se llega a la conclusi\'on (\emph{horrorosa, pero incuestionable}) de que existe al menos un conjunto infinito que no contiene un subconjunto numerable. Sin embargo tambi\'en existe la demostraci\'on de lo contrario. Pero, como la demostraci\'on de que, en efecto, todo conjunto infinito posee un subconjunto numerable es indiscutible, se produce la tendencia natural a cuestionar la demostraci\'on que aqu\'i se ha presentado.\\ 
En ese sentido hay dos elementos que deben valorarse con especial atenci\'on. En caso de haber un error en la demostraci\'on del teorema $1$, debe sospecharse, en primer lugar, de la suposici\'on, realizada en ella, de la existencia de alg\'un conjunto infinito que no es $p-conjunto$ lo cual coincide, una vez hecho el cambio de denominaci\'on, con la suposici\'on de que existe un conjunto infinito que no es numerable y, en segundo lugar, el paso al l\'imite, tambi\'en efectuado en esta demostraci\'on, como procedimiento que puede darse por concluido y que permite la ``extracci\'on'' de \textbf{\emph{TODOS}} los, determinados, subconjuntos de cierto conjunto.\\
En definitiva, para escapar de la paradoja ``Todo conjunto infinito tiene y no tiene un subconjunto numerable'' debe aceptarse una proposici\'on evidente: ``Todos los conjuntos infinitos son numerables'', o por el contrario es il\'icito en general el paso al l\'imite en teor\'ia de conjuntos.\\
Aqu\'i, la idea de que todos los conjuntos infinitos sean numerables es al parecer la m\'as d\'ebil; sin embargo no debe desech\'arsele tan a la ligera. T\'engase en cuenta que, como alternativa, se tiene una aceptaci\'on expresa del rechazo a considerar a los procesos que involucran al infinito como algo acabado en alg\'un momento, y la desconfianza consecuente respecto a los pasos al l\'imite tan necesarios y generalizados. En tal caso, por ejemplo, una vez dado $A=(\bigcup_{i=1}^\infty Pc_i )\cup(lim_{n\longrightarrow\infty}A_n )$, el conjunto $lim_{n\longrightarrow\infty}A_n $ seguir\'a conteniendo al menos un $Pc$ y pasando al l\'imite, cuantas veces se quiera, seguir\'a sucediendo lo mismo, por lo que nos enfrentamos a un procedimiento nunca acabado que contradice la idea del infinito actual, por lo que, de alguna manera, si se acepta la idea del infinito actual se acepta al mismo tiempo que todos los conjuntos infinitos son numerables.\\ 
El infinito actual suele ser, sin embargo, el reflejo, en las matem\'aticas, de la verdad indiscutible que expresa la infinitud del mundo en general, infinitud que no puede ser concebida por mentes estrechas, que no pueden imaginar una verdad m\'as all\'a de la que pueden construir con sus propios pensamientos.\\
Del teorema $2$ se deduce que todas las potencias de $\mathbb{N}$ son numerables y se produce el dilema siguiente: $\mathbb{R}$ es numerable, o no es la potencia del conjunto $\mathbb{N}$. Es decir ser\'a $c=\aleph_0$ o por el contrario $c\neq 2^{\aleph_0}$. Sin embargo se ha realizado la propuesta del conteo de los reales de $\left[0,1\right]$ ``por capas''. De aqu\'i que todas las potencias de $\mathbb{R}$ sean numerables.\\
N\'otese, por otra parte, que la situaci\'on contradictoria a la que se llega en el caso del teorema de la imposibilidad de establecer una relaci\'on $uno-a-uno$ entre los elementos de un conjunto dado y sus subconjuntos (Cuarto caso), es totalmente equivalente a lo que sucede con la paradoja del barbero y todas las paradojas an\'alogas a ella. Por lo que se puede afirmar que el problema que introdujo a las matem\'aticas en una gran crisis, es simplemente el resultado del error. Por lo que habr\'ia que ver hasta qu\'e punto han sido necesarios todos lo esfuerzos para deshacerse de las paradojas.\\
Por \'ultimo quiero referirme, de forma somera, al problema de ``la medida''. Es indudable que en base a los resultados que he obtenido aqu\'i, todos los conjuntos num\'ericos y sus productos cartesianos ser\'an de medida cero. Este problema, sin embargo, ser\'a tratado en posteriores  investigaciones y para nada es indicio de contradicci\'on en los an\'alisis aqu\'i efectuados que se basan en demostraciones evidentes. Es m\'as razonable, en vista de lo que se ha expuesto, que el c\'umulo de errores l\'ogicos sea bastante grande en la teor\'ia actual que tiene como base la aceptaci\'on de absurdos como los tratados en los casos dese el uno hasta el cuatro de la segunda secci\'on de este art\'iculo. 

\end{document}